\documentclass[11pt,letterpaper,reqno]{amsart}

\usepackage[final
]{showkeys}

\allowdisplaybreaks

\usepackage{AKstyle}
\usepackage[toc,page]{appendix}
\numberwithin{equation}{section}

\usepackage{amsmath}
\usepackage{caption}
\usepackage[labelfont=rm]{subcaption}
\usepackage{scalerel,stackengine}
\stackMath
\newcommand\reallywidehat[1]{%
\savestack{\tmpbox}{\stretchto{%
 \scaleto{%
 \scalerel*[\widthof{\ensuremath{#1}}]{\kern-.6pt\bigwedge\kern-.6pt}%
 {\rule[-\textheight/2]{1ex}{\textheight}}
 }{\textheight}%
}{.8ex}}%
\stackon[1pt]{#1}{\tmpbox}%
}
\usepackage[pagebackref=true, colorlinks]{hyperref}

\renewcommand{\phi}{\varphi}

\newcommand{\mat}[1]{\left(\begin{matrix} #1 \end{matrix} \right)}  

\hypersetup{pdffitwindow=true,linkcolor=blue,citecolor=blue,urlcolor=blue,menucolor=blue}

\usepackage{comment}

\begin{document}

\author{Xin Zhang}
\thanks{The author is supported by GRF grant 17300324 and NSFC grant 12001457.}
\email{xz27@hku.hk}
\address{Department of Mathematics, The University of Hong Kong}

\title[Zaremba's Conjecture]
{Expansion in $\text{SL}_2(\mathbb Z/q\mathbb Z)$ and Zaremba's Conjecture}

\begin{abstract} We establish an expansion theory for $\text{SL}_2(\mathbb Z/q\mathbb Z)$. Incorporating this into a framework recently developed by Shkredov, we confirm Zaremba's conjecture. 
\end{abstract}
\date{\today}
\maketitle

\section{Introduction}

Every rational number $\frac{a}{q}\in (0,1)$ admits a continued fraction expansion: 
$$x=\frac{1}{a_1+\frac{1}{a_2+\ddots \substack{\\\\+\frac{1}{a_s}}}}:=[a_1,a_2,\cdots, a_s],$$
with $a_1,a_2,\cdots, a_s\geq 1$. The integers $a_i$ are called partial quotients of $a/q$. \par

In 1971, Zaremba formulated the following assertion: 
\begin{conj}[\normalfont Zaremba, 1972 \cite{Zar72}, p. 76]\label{Zaremba} Every natural number is the denominator of a reduced fraction whose partial quotients are absolutely bounded. That is, there exists some $\mathcal M>0$ such that for any $q\in\mathbb N$, one can find a natural number $a<q$ with $\normalfont\text{gcd}(a,q)=1$, such that writing $a/q=[a_1,a_2,\cdots, a_s]$, we have $a_i\leq \mathcal M$ for every $1\leq i\leq s$. 
\end{conj}

In this paper we prove 
\begin{theorem}\label{main}
Zaremba's conjecture holds. 
\end{theorem} 

In fact, Zaremba conjectured that $\mathcal M=5$. Our method can only produce a large but an absolute bound.\par 
In a spectacular paper \cite{BK14a}, Bourgain-Kontorovich proved that $\mathcal M=50$ works for \emph{almost} all natural numbers. Later, the bound $\mathcal M=50$ in this density one theorem was improved to $\mathcal M=5$ by Huang \cite{Huang15}. For $\mathcal M=4$, Kan obtained a positive density result \cite{Kan15}.   In another direction, Korobov in 1960 proved that for all natural numbers $q$, there exists $a<q$ such that all partial quotients of $\frac{a}{q}$ are bounded by $O(\log q)$ \cite{Kor60}. This bound was improved to $O(\frac{\log q}{\log\log q})$ for a specific subset of natural numbers by Murphy-Moshchevitin-Shkredov \cite{MMS26}, and in particular, this set contains all primes. In a recent brilliant paper by Shkredov \cite{Sh26}, this bound was further improved to an absolute constant $\mathcal M$.  Both \cite{MMS26} and \cite{Sh26} are considering some cantor type sets with some additive structures in $\mathbb Z/q\mathbb Z$ that are originally due to Moshchevitin \cite{Mo07}, and a critical ingredient in these two works is Helfgott's triple product theorem on $\text{SL}_2(\mathbb Z/p\mathbb Z)$ \cite{He08} as well as its generalization to square free moduli by \cite{BGS10}. \par

Nevertheless, the set of integers covered by \cite{MMS26} and \cite{Sh26} still only makes a positive proportion of all natural numbers. As hinted in the concluding remarks of \cite{MMS26}, to deal with an arbitrary modulus $q$, one naturally needs a growth theory in $\text{SL}_2(\mathbb Z/q\mathbb Z)$ for a general $q$. In \S 2, we provide such a theory (Theorem \ref{product}, Theorem \ref{bd}, Theorem \ref{measure}), albeit we have to admit it is not light. In \S 3 we apply the results in \S 2 to estimate the number of solutions to a particular equation in $\mathbb Z/q\mathbb Z$ to obtain Theorem \ref{count}. Implementing Theorem \ref{count} in the framework provided by the proof of Theorem 8 in \cite{Sh26} then leads to Theorem \ref{main}. For readers' convenience, in \S \ref{shk} we give an outline of Shkredov's argument how to use Theorem \ref{count} to deduce Theorem \ref{main}. We invite the readers to Shkredov's work \cite{Sh26} for details.

\begin{remark}
The bound $\mathcal M$ in Conjecture \ref{Zaremba} is effectively computable. The methods in \cite{Sh26} and in this paper would produce a large number. Lowering this bound to 5 or any reasonable magnitude requires a new idea. 
 \end{remark}
 
\begin{remark} Bourgain-Varj\'u proved a spectral gap theorem in $\text{SL}_d(\mathbb Z/q\mathbb Z)$ with respect to a fixed finite generating set $S$ and an arbitrary modulus $q$. However, in the consideration of Zaremba's conjecture, we need to deal with a generating set $S$ that depends on $q$ (see \eqref{0452}), in particular, $|S|\approx q^{2\tau}$ for some small but absolute $\tau>0$, and the dependence of the spectral gap produced in \cite{BV12} is not good enough for our consideration here. We naturally need a theory in the spirit of Helfgott's triple product theorem in $\text{SL}_2(\mathbb Z/p\mathbb Z)$ \cite{He08}.
\end{remark}

\noindent {\bf Acknowledgement} The author is grateful to Alex Kontorovich for discussions on this work.

\section{Expansion in $\text{SL}_2(\mathbb Z/q\mathbb Z)$, $q$ arbitrary}

In this section we give some theorems on the expansion of $\text{SL}_2(\mathbb Z/q\mathbb Z)$ for a general modulus $q$. For notational convenience, sometimes we write $\Gamma=\text{SL}_2(\mathbb Z), \Gamma(q)$ is the principal congruence subgroup of $\Gamma$ of level $q$, and $\Gamma_q=\Gamma/\Gamma(q)$. \par
The first theorem is on the triple product growth of sets in $\text{SL}_2(\mathbb Z/q\mathbb Z)$. The special case $q$ square free was proved as Proposition 4.3 in \cite{BGS10}.
\begin{thm}[triple product growth]\label{product} Fix $\kappa_1, \kappa_2>0$. There exists $\omega, \delta>0$ depending on $\kappa_1, \kappa_2$ such that the following holds: Let $q\in \mathbb Z_+$ and $A$ be a subset of $\normalfont\text{SL}_2(\mathbb Z/q\mathbb Z)$ satisfying the following properties:
\begin{enumerate}
\item $q^{\kappa_1}<|A|<q^{3-\kappa_1};$
\item For all $t\in \mathbb Z/q\mathbb Z$ and for all $g\in\normalfont{\text{Mat}}_2(q)$ with $\pi_p(g)\neq \mat{0&0\\0&0}$ for all primes $p|q$, we have for all $q_1| q, q_1>q^{\omega}$,
$$\#\{x\in A|  \normalfont{\text{Tr}}(gx)-t=0(q_1) \}<q_1^{-\kappa_2}(|A|).$$
\end{enumerate}
Then 
$$|A\cdot A\cdot A|>q^{\delta}|A|.$$
\end{thm} 
Roughly speaking, Theorem \ref{product} says that if a set $A\subset \text{SL}_2(\mathbb Z)$ is not significantly concentrated in any small linear subvariety, and the size of $A$ is not too large, then we expect a significant growth from $|A|$ to $|A\cdot A\cdot A|$. \par
Closely related to Theorem \ref{product} are the following two theorems regarding bounded generation and $l^2$-flattening:

\begin{thm}[bounded generation]\label{bd} Fix $\kappa >0$. For any $\epsilon >0$, there exists $\omega>0, K\in\mathbb Z_+$ depending on $\kappa, \epsilon$ such that the following holds: \par
Let $q\in \mathbb Z_+$ and $A$ be a subset of $\normalfont\text{SL}_2(\mathbb Z/q\mathbb Z)$ satisfying: For all $t\in \mathbb Z/q\mathbb Z$ and for all $g\in\normalfont{\text{Mat}}_2(q)$ with $\pi_p(g)\neq \mat{0&0\\0&0}$ for all primes $p|q$, we have for all $q_1| q, q_1>q^{\omega}$,
\begin{align}\label{0722}
\#\{x\in A|  \normalfont{\text{Tr}(gx)}-t=0(\mod q_1)\}<q^{-\kappa}(|A|).
\end{align}
Then 
$$A^{K}\supset \Gamma(Q)/\Gamma(q),$$
where $Q|q$ and $Q<q^\epsilon$.
\end{thm} 

\begin{thm}[$l^2$-flattening]\label{measure}
Given $\gamma, \kappa >0$, there exists $\omega, \eta>0$ such that the following holds: Let $q\in\mathbb Z_+$. Let $\mu\in {\normalfont \text{Prob}}(\Gamma_q)$ such that $$\Vert \mu\Vert_2>q^{-\frac{3}{2}+\gamma},$$
and for all $t\in \mathbb Z/q\mathbb Z$ and for all $g\in\normalfont{\text{Mat}}_2(q)$ with $\pi_p(g)\neq \mat{0&0\\0&0}$ for all primes $p|q$, we have for all $q_1| q, q_1>q^{\omega}$,
\begin{align}\label{0426}
\mu (\{x\in A| \normalfont{\text{Tr}}(gx)-t=0(\mod q_1) \})<q_1^{-\kappa}.
\end{align}
Then $$\Vert \mu *\mu\Vert_2<q^{-\eta} \Vert \mu\Vert_2.$$

\end{thm}

Theorem \ref{bd} can easily imply Theorem \ref{product} through the Helfgott inequality (Lemma 2.2, \cite{He08}): For any $l\geq 3$, $$\frac{|A|^l}{|A|}\leq \left (\frac{|A^3|}{|A|}\right)^{l-2}.$$ Also, Theorem \ref{product} implies Theorem \ref{measure}, following the argument in the proof of Proposition 4.2, \cite{BGS10}, which proved the case that $q$ is square free for Theorem \ref{measure} but the proof works for a general modulus $q$. So we focus on proving Theorem \ref{bd}. 

\subsection{Square Free Case} We first prove Theorem \ref{bd} for $q$ square free. Theorem \ref{product} and Theorem \ref{measure} for $q$ square free are proved in the paper \cite{BGS10}.\par
 Let $A$ be as given in Theorem \ref{bd}. Set $\mu=\mu_A$, where $\mu_A$ is the uniform counting measure supported on $A$. \par 
We apply Theorem \ref{measure} to $\mu$. Set $\gamma=\frac{\epsilon}{40}$, and let $\omega, \eta$ be the implied constants by $\kappa, \gamma$ in Theorem \ref{measure}. If $\Vert \mu\Vert_2\leq q^{-\frac{3}{2}+\gamma}$, then by Cauchy-Schwarz we obtain $|A|=|\text{Supp}(\mu)|\geq q^{3-2\gamma}$. If $\Vert \mu\Vert_2> q^{-\frac{3}{2}+\gamma}$, we can apply Theorem \ref{measure} to obtain $\Vert \mu*\mu\Vert_2\leq  q^{-\eta}  \Vert \mu\Vert_2$. Observe that the non-concentration condition \eqref{0426} is preserved if we replace $\mu$ by $\mu*\mu$, with same parameter $\omega, \kappa$. Iteratively taking convolutions of $\mu$ until we obtain a number $l=O(\frac{1}{\eta})$ such that $\Vert \mu^{(2^l)} \Vert_2\leq q^{-\frac{3}{2}+\gamma}$. Thus by Cauchy-Schwarz, $|A^{2^l}|=|\normalfont\text{Supp}(\mu^{(2^l)})|>q^{3-2\gamma}$. An application of Proposition 3.19, \cite{TZ23b} to $A^{2^l}$ then yields Theorem \ref{bd}, for $q$ square free. \par

\subsection{Bounded Exponent Case} In this section we prove Theorem \ref{bd} for the case $q=\prod_{i=1}^m p_i^{n_i}$, $n_i\leq L, \forall i$. \par

Let $A, \kappa, \epsilon$ be as given in Theorem \ref{bd}. We take $\omega=\frac{\omega_0(\kappa, \frac{\epsilon}{80L})}{L}$, where $\omega_0$ is the implied function for $\omega$ in Theorem \ref{measure}. 

Let $q_0$ be the square free part of $q$. We take $\mu=\pi_{q_0}^*(\mu_A)\in \text{Prob}(\Gamma_{q_0}) $, the pushforward of $\mu_A$ under the map $\pi_{q_0}: \Gamma_q\rightarrow \Gamma_{q_0}$.  Then $\mu$ satisfies the $(\omega_0, \kappa)$ non-concentration property. An application of Theorem \ref{measure} gives $\Vert \mu^{(2^{l})}\Vert_2 \leq q_0^{-\frac{3}{2}+\frac{\epsilon}{80L}}$, $l=O_{\kappa, \epsilon}(1)$, which implies  $|\pi_{q_0}(A^{2^l})|>q^{3-\frac{\epsilon}{40L}}$. It then follows from Proposition 3.19 \cite{TZ23b} that there exists $Q_0<q_0^{\frac{\epsilon}{2L}}$, such that 
 
 \begin{align*}
 \pi_{q_0}(A^{O_{\kappa,\epsilon}(1)  })\supset \Gamma(Q_0)/\Gamma(q_0).
  \end{align*}
 Then we can follow the argument in the proof of Proposition 3, \cite{BV12} to show that
 $$\pi_{q_0}(A^{O_{L, \kappa,\epsilon}(1)  })\supset \Gamma(Q_1)/\Gamma(q),$$
for some $Q_1|q, Q_1<q^\epsilon$ as required.

\subsection{Large Exponent Case} In this section we consider $q=\prod_{i=1}^m p_i^{n_i}$, where $\forall i, n_i>L$ for a sufficiently large $L$. 

\begin{proposition}\label{key} Let $\kappa, \epsilon$ be as given in Theorem \ref{bd}.  There exists $L=O(\frac{1}{\epsilon})\in \mathbb Z_+$, $\omega=\omega(\kappa, \epsilon)>0$, $C=C(\kappa, \epsilon)\in\mathbb Z_+,$ such that for any $q=\prod_{i=1}^m p_i^{n_i}, n_i>L, \forall i$, any $A$ satisfies \eqref{0722}, we have
$$A^{C}\supset \Gamma(q_*^{\{\epsilon\}})/\Gamma(q_*),$$
where $q_*\Vert q$ and $q_*\geq q^c$ for some $c>0$ depending only on $\kappa$. 
\end{proposition}

Here $q^*\Vert q$ means for any $p|q^*$, $v_p(q^*)=v_p(q)$, where $v_p$ is the valuation function. If $q=\prod_{i=1}^m p_i^{n_i}$, the notation $q^{\{\epsilon\}}$ means $\prod_{i=1}^m p_i^{[n_i\epsilon]}$. 

In our view Proposition \ref{key} may be a most technical part of this expansion theory. It relies on a sum-product theorem over $\mathcal O/\mathcal I$, where $\mathcal O$ is the ring of integers of a general number field $\mathcal K$ and $\mathcal I$ is a general ideal of $\mathcal O$. This sum product theorem was first conjectured in \cite{SG20} and was later proved as Theorem 1.3 in \cite{TZ25}. \par

To prove Proposition \ref{key}, we apply the special case of the aforementioned sum-product theorem where $\mathcal K$ is a quadratic extension of $\mathbb Q$, to the framework provided by \cite{BG09} to exploit the Lie structure $\mathfrak{sl}_2(\mathbb Z/q\mathbb Z)$. The proof of Proposition \ref{key} is virtually identical to that of Proposition 5.1, \cite{TZ23b}, and we invite the readers to consult \cite{TZ23b} for details. \par

\subsection{General Case} For a general $q$, write $q=\prod_{i=1}^m p_i^{n_i}=q_sq_l$, where $q_s=\prod_{i: n_i\leq L}p_i^{n_i}$ and $q_l=\prod_{i: n_i> L}p_i^{n_i}$. Let us assume that both $q_s$ and $q_l$ are not too small. Theorem \ref{bd} and Proposition \ref{key} only tells us that under the reductions of certain moduli, a bounded product of $A$ has almost full dimension. It naturally raises the question: Given $q_1, q_2\Vert q, (q_1, q_2)=1$ and a set $A\subset \Gamma_q$ such that $\pi_{q_1}(A), \pi_{q_2}(A)$ are large, can one show that a bounded product of $A$ is large under the projection of $q_1q_2$? We do have such a gluing tool to help achieve this, which is precisely stated as follows: 
\begin{proposition}\label{glu}Fix $0<\theta<\min\{10^{-12}, \frac{1}{L}\}$. Suppose $q_1\Vert q$, $q_2\Vert q_l$, and $(q_1, q_2)=1$. Suppose for some set $A\subset \Gamma_q$ we have 
$|\pi_{q_1}(A)|>q_1^{3-\theta}, |\pi_{q_2}(A)|>q_2^{3-\theta}$. Then there exists $q_3\Vert q_2$, $q_3>q_2^{\frac{1}{4}10^{-4}}$, such that 
$$|\pi_{q_1q_3} (A^{O(\log \frac{1}{\theta})^2} ) |>(q_1q_3)^{3-O(\theta^{\frac{1}{2}})}.$$

\end{proposition}
For a proof of Proposition \ref{glu}, we refer the reader to Proposition 6.1, \cite{TZ23b}, which deals with a more complicated group $\text{SL}_2(\mathbb Z/q\mathbb Z)\times \text{SL}_2(\mathbb Z/q\mathbb Z)$. Since any homomorphism $f: \Gamma_{p_i^{n_i}}\rightarrow \Gamma(p_j)/\Gamma(p_j^{n_j})$ is trivial if $p_i\neq p_j$, $p_j\neq 2$, and is trivial mod $p_j^{n_j-1}$ if $p_i\neq p_j$, $p_j= 2$, compared to \cite{TZ23b}, we only need to deal with two cases here, and the arguments in \S 6.1 and \S 6.2 in \cite{TZ23b} dealing with Event 1 and Event 2.1 suffice to give a proof for Proposition \ref{glu}. This also explains a change of exponent in $O(\theta^{\frac{1}{2}})$ compared to  $O(\theta^{\frac{1}{4}})$ in Proposition 6.1 of \cite{TZ23b}. \par
Theorem \ref{bd} can then be obtained by iteratively applying Proposition \ref{key} and Proposition \ref{glu}. Indeed, observing that the growth rate of modulus in each iterative step, described by $c$ in Proposition \ref{key} and $\frac{1}{4}10^{-4}$ in Proposition \ref{glu} are some positive constants depending only on the non-concentration parameter $\kappa$, it ensures that a bounded product of $A$ can cover a set of almost full density in $\Gamma_q$. See \S 7 of \cite{TZ23b} for details. 

\section{Counting solutions to an equation}

Our main goal in this section is Theorem \ref{count}. 

As noted in the paper \cite{Sh26}, Zaremba's conjecture is related to counting solutions to the following equation 
\begin{align}\label{0102}
(a+2j)(b+2j)=1(\mod q),
\end{align}
where $a\in A, b\in B$, and $A, B$ are Cantor type sets of the form 
\begin{align}\label{1002}
A=\Lambda_1\dotplus [1, N], B=\Lambda_2\dotplus [1, N]
\end{align}
 in $\mathbb Z/q\mathbb Z$. Here for two sets $S_1, S_2$, $S_1\dotplus S_2$ means for $x_1, y_1\in S_1, x_2, y_2\in S_2$, $x_1+x_2=y_1+y_2 \Rightarrow x_1=x_2, y_1=y_2$. The sets $A, B$ roughly have size $t^{2w_M} N$, where $ t^2 N=q, t\sim q^{\frac{1}{2}-\tau}, N\sim q^{2\tau}$, and $w_M$ is the Hausdorff dimension of the set of reals in $[0,1]$ with all partial quotients bounded by $M$. By a theorem of Hensley \cite{Hen92}, $1-w_M\sim \frac{1}{M}$ as $M\rightarrow \infty$.  \par

The equation \eqref{0102} is equivalent to 
\begin{align}\label{0955}
\mat{-2j& 1-4j^2\\1& 2j}a=b,
\end{align}
where we are considering the linear fractional transformation of $\Gamma_q$ on $\mathbb P^1(\mathbb Z/q\mathbb Z)$, which is transitive.  \par
We prefer to work with $\text{SL}_2(\mathbb Z)$, so we consider the following equation which is equivalent to $\eqref{0955}$:
\begin{align}\label{0103}
g_j a =-b,
\end{align}
where $g_j=\mat{2j& 1-4j^2\\-1& 2j}$. \par 
Let 
\begin{align}\label{0452}
S=\{g_j: j=1, \cdots ,N\}
\end{align}
Since $-B$ has a similar structure as $B$, we aim to prove: 
\begin{theorem}\label{count}  Let $A, B, S$ be as in \eqref{1002} and \eqref{0452}. We have  
\begin{align}\label{0905}
\sum_{s\in S}\sum_{a\in A}\sum_{b\in B}{\bf 1}\{sa=b\}= \frac{ \phi(q)}{q^2}|A||B|N+O_\epsilon(q^{\epsilon} \sqrt{|A||B|}N^{1-\eta}),
\end{align}
where $\eta>0$ is some absolute constant and $\epsilon$ is arbitrarily small. 

 \end{theorem}
 
 In particular, the error term in \eqref{0905} is smaller than the main term if we take 
 
 \begin{align}\label{1007}
 \tau \eta> \frac{1}{M}. 
 \end{align}
 \par 
The inequality \eqref{1007} puts a main restriction on $M$. The constant $\eta$ we produce is very small, which turns into a large lower bound for $M$. 
 
Compare Theorem \ref{count} with Lemma 13 in \cite{Sh26}. While we agree with Lemma 13 when $q$ is a prime, our main term is different from the main term $\frac{|A||B|}{q}$ if $q$ is a composite number with small prime divisors.  \par
Let $\mu=\mu_S$ be the normalized counting function on $\Gamma_q$. \par
To prove Theorem \ref{count}, we can do a spectral decomposition for $\mu$, following the idea in \cite{MMS26}. However, unlike the case that $q$ is a prime, if $q$ is a general modulus, we need to have a more careful analysis on non-constant low-frequency harmonics, and it would be too wasteful to throw them away to the error term. \par
For $Q|q$, define $\mu_Q(x): \Gamma_q\rightarrow \mathbb R_{\geq 0}$ as $$\mu_Q(x) =\mathbb E_{y\equiv x(\mod Q)} (\mu(y))=\frac{1}{|\Gamma(Q)/\Gamma(q)|}\sum_{\substack{y\in \Gamma_q\\ y\equiv x (\mod Q) }}\mu(y).$$
So $\mu_Q$ is constant on each congruence class mod $Q$. \par

Write $Q=\prod_{i=1}^t p_i^{n_i}$, and consider the expression 
$$\prod_{i=1}^t(p_i^{n_i}-p_{i}^{n_i-1})=\sum_{q' | Q} m(q') q', $$ where $m(q')=0, \pm 1$. \par
Let $$f_Q(x)=\sum_{q'|Q}m(q')\mu_{q'}(x).$$
For example, $f_2= \mu_2-\mu_1$ and $f_{12}=\mu_{12}-\mu_6-\mu_4+\mu_2$. \par
It is an exercise to check that 
\begin{align}\label{0327}
\mu(x)=\sum_{Q|q}f_Q(x),
\end{align}
and each $f_Q$ lives in the space $H_Q\subset l^2(\Gamma_q)$, where $H_Q$ consists of all functions $g$ such that $g$ is constant over any congruence class mod $Q$, and for any $q'|Q, q'\neq Q$, the summation of $g$ over any congruence class mod $q'$ is 0. It is clear that $H_Q$ is a subrepresentation of the right regular representation of $\Gamma_q$ on $l^2(\Gamma_q)$ and is not induced by any smaller divisor of $Q$. \par
Now for each $Q$ small, namely, $Q<N^{\frac{1}{2}}$, write $q=q_1q_2$, where $(q_1, q_2)=1$ and $q_1, Q_1$ have the same set of prime divisors. \par
We can estimate 
\begin{align}
&\sum_{s\in \Gamma_q}\mu_Q(s) \sum_{x\in A} B(sx)\nonumber\nonumber \\ 
\label{0243}=&\sum_{\bar a, \bar b (\mod Q)} \sum_{\bar {s}\in \Gamma_Q}{\bf 1}\{\bar{s} \bar{a}= \bar{b}\} \sum_{\substack{a\in A, b\in B\\ a\equiv \bar a (\mod Q) \\ b\equiv \bar b(\mod Q)}} \sum_{\substack {s\in \Gamma_q\\ s\equiv \bar s (\mod Q)}}\mu _Q(s){\bf 1 }\{sa=b\}
\end{align}

Write $\bar s=\mat{x&y\\z&w}$. The equation $\bar s \bar a=\bar b$ can then be written as $x\bar a+y-z\bar b\bar a-w\bar b=0$. \par
We observe that the Jacobian
$$\frac{\partial (x\bar a+y-z \bar a\bar b-w\bar b, x w-yz)}{\partial (x, y, z, w)}= \mat{\bar a&1&-\bar a\bar b&-\bar b\\ w&-z&-y&x} $$
has rank 2 $(\mod p)$ for each $p|Q$, as $\bar s \bar a=\bar b$ implies $ z\bar a+w\neq 0(\mod p)$ for each $p|Q$. Also $\mu_Q(s)$ is constant in the last sum of \eqref{0243}. Write this constant as $\mu_Q(\bar{s})$. Hensel's lifting then gives 
\begin{align}\label{0255}
\eqref{0243}=\sum_{\bar{a}, \bar{b}(\mod Q)}\sum_{\bar s \in \Gamma_Q}{\bf 1}\{\bar{s} \bar{a}= \bar{b}\} \left( \frac{|A||B|}{Q^2}(1+O(\frac{1}{N})) \mu_Q(\bar s)(\frac{q_1}{Q})^2\frac{|\Gamma_{q_2}|}{|\mathbb P^1(\mathbb Z/q_2\mathbb Z)|}  \right)
\end{align}
Observe that for each $\bar{s}$ where $\mu_Q(\bar{s})\neq 0$, there are exactly $\phi(Q)$ many $\bar a$ such that $\bar s\bar a\in \mathbb Z/Q\mathbb Z$ (see the definition of $S$ at \eqref{0452}). So when $\bar{s}, \bar a$ are fixed, $\bar b$ is then determined, so 
\begin{align}\label{0319}
\eqref{0255}=&\sum_{\bar s \in \Gamma_Q}\mu _Q(\bar s)  \left( \frac{\phi(Q)|A||B|}{Q^2}(1+O(\frac{1}{N})) (\frac{q_1}{Q})^2\frac{|\Gamma_{q_2}|}{|\mathbb P^1(\mathbb Z/q_2\mathbb Z)|}  \right) \nonumber \\
= &\frac{\phi(Q)|A||B|}{|\Gamma(Q)/\Gamma (q)| Q^2}(1+O(\frac{1}{N})) (\frac{q_1}{Q})^2\frac{|\Gamma_{q_2}|}{|\mathbb P^1(\mathbb Z/q_2\mathbb Z)|}  \nonumber\\
= &\frac{|A||B|}{|\mathbb P^1(\mathbb Z/q\mathbb Z)|}\prod_{p|Q}(1-\frac{1}{p^2})(1+O(\frac{1}{N}))
\end{align}

We can then deduce from \eqref{0319} that
\begin{align}
\sum_{s\in \Gamma_q}f_Q(s)\sum_{x\in A}B(sx)=\nu_Q\frac{|A||B|}{|\mathbb P^1(\mathbb Z/q\mathbb Z)|}+ O(\frac{|A||B|}{|\mathbb P^1(\mathbb Z/q\mathbb Z)|N^{1-\epsilon}}),
\end{align}
where $\nu_Q=\frac{ (-1)^{\omega (Q)}}{Q^2}$ if $Q$ is square free and 0 otherwise. \par
In view of the decomposition \eqref{0327}, we obtain 
\begin{align}\sum_{s\in \Gamma_q}\mu (s)\sum_{x\in A}B(sx) = \mathfrak M+\mathfrak E, \end{align}
where 
\begin{align}\nonumber \mathfrak M:=& \sum_{s\in \Gamma_q}\sum_{\substack {Q|q\\Q<N^{\frac{1}{2}}}}f_Q(s)\sum_{x\in A}B(sx)\\ \nonumber 
= &(\sum_{\substack{Q|q\\Q<N^{\frac{1}{2}}}}\nu_Q)\frac{|A||B|}{|\mathbb P^1(\mathbb Z/q\mathbb Z)|} +O(\frac{|A||B|}{|\mathbb P^1(\mathbb Z/q\mathbb Z)|N^{1-\epsilon}}). \\
=&\frac{\phi(q)}{q^2}{|A||B|}+O(\frac{|A||B|}{qN^{1-\epsilon}}) \label{0909}
\end{align}
and 
$$\mathfrak E:= \sum_{\substack {Q|q\\Q\geq N^{\frac{1}{2}}}}\sum_{s\in \Gamma_q} f_Q(s)\sum_{x\in A}B(sx)$$
For each $Q\geq N^{\frac{1}{2}}$, we have 
\begin{align} \label{0837}
\vert \sum_{s\in \Gamma_q}\sum_{x\in A}f_Q(s) B(sx) \vert \leq \sqrt{|A||B|}\left(|A|^{-1}\sum_{s}r_{f_Q, 2^l}(s)\sum_{x\in A}A(sx)\right)^{\frac{1}{2^l}},
\end{align}
where $$r_{f_Q, 2^l}(s)=\sum_{\substack {s_1, \cdots, s_{2^l} \in \Gamma_Q \\ s_1s_2^{-1}\cdots s_{2^l-1}s_{2^l}^{-1}=s }}  f_{Q}(s_1)\cdots f_Q(s_{2^l}).$$ 
This follows from iterative applications of the Cauchy-Schwarz inequality. See Lemma 9, \cite{MMS26}. The point is, we can write 
\begin{align}\label{0950}
r_{f_Q, 2^l}= (TT^*)^{2^{l-1}}  (f_Q*\hat{f_Q}), 
\end{align}
 where $\hat{f_Q}(s)=f_Q(s^{-1})$, and the two operators $T, T^*: l^2(\Gamma_q)\rightarrow l^2(\Gamma_q)$ are defined as: for $g\in l^2(\Gamma_q)$, 
 \begin{align*}
 &T(g)(x)=\sum_{s\in G}\mu_Q(s)g(sx),\\
 &T^*(g)(x)=\sum_{s\in G}\hat{\mu_Q}(s)g(sx),
 \end{align*}
 where $\hat{\mu_Q}(s)=\mu_Q(s^{-1})$. \par
 The equation \eqref{0950} follows from a simple calculation once observing that $$\mu_{q_1}*\mu_{q_2}=\mu_{\text{gcd}(q_1, q_2)}*\mu_{\text{gcd}(q_1, q_2)}.$$ W can thus write $(TT^*)^{(n)}(g)(x)=\sum_{s\in \Gamma_Q} (\mu_Q*\hat{\mu_Q})^{(n)}(s) g(sx ) $.  \par

We first take an integer $k_0=O(\frac{\log Q}{\log N})$ if $N<Q$, and $k_0=1$ if $N\geq Q$, so that $(\mu_Q*\hat{\mu_Q})^{(k_0)}$ satisfies the non-concentration property \eqref{0426} described by a constant $\kappa=\kappa_0>0$, which is absolute. For the existence of such an absolute constant $\kappa_0$, see Page 10-11 of \cite{MMS26}. 

By $l^2$-flattening (Proposition \ref{measure}), there exists $k_1=O(1)$ such that $$ \Vert (\mu_Q*\hat{\mu_Q})^{(k_0k_1)} \Vert_2 \leq \frac{Q^{\frac{1}{4}}}{|\Gamma_q|}.$$ 
Let $d(Q)$ be the minimal dimension of a faithful irreducible subrepresentation of the right regular representation of $\Gamma_Q$. By Lemma 7.1, \cite{BG08a}, we have $d(Q)\geq\frac{Q}{3}$, so by considering the trace of $(TT^*)^{k_0k_1}$, we obtain
\begin{align*}
&\Vert (\mu_Q*\hat{\mu_Q})^{(k_0k_1)} *(f_Q*\hat{f_Q}) \Vert_2 \\\leq & \frac{ \Vert (\mu_Q*\hat{\mu_Q})^{(k_0k_1)} \Vert_2 \sqrt{|\Gamma_q|}}{\sqrt{d(Q)}}\Vert f_Q*\hat{f_Q}\Vert_2\\\leq& {Q^{-\frac{1}{4}} }   \Vert f_Q*\hat{f_Q}\Vert_2   .
\end{align*}
Iteratively applying the operator $(TT^*)^{k_0k_1}$ to $f_Q*\hat{f_Q} $ for $k_2$ times for some $k_2\asymp \frac{\log q}{\log Q}$, we obtain 
$$\Vert (\mu_Q*\hat{\mu_Q})^{(k_0k_1k_2)} *(f_Q*\hat{f_Q}) \Vert_2\ll \frac{1}{q^8},$$
which surely implies 
$$ (\mu_Q*\hat{\mu_Q})^{(k_0k_1k_2)} *(f_Q*\hat{f_Q})(s)\leq \frac{1}{q^4} $$
for every $s\in \Gamma_q$. \par

In \eqref{0837} we can take $l$ so that $2^l\approx k_0k_1k_2$ such that the bound $r_{f_Q, 2^l}(s)\leq \frac{1}{q^4}$ holds. We then obtain 
$$\text{RHS of }\eqref{0837} \leq \sqrt{|A||B|}N^{-\eta}$$
for some absolute $\eta>0$.  Therefore, 
\begin{align}\label{0908}
\mathfrak E\ll_\epsilon q^{\epsilon} \sqrt{|A||B|}N^{-\eta}.
\end{align}
for an arbitrary $\epsilon>0$. \par
Theorem \ref{count} is thus proved by combining \eqref{0909} and \eqref{0908}. 

\section{An outline of Shkredov's argument}\label{shk}
In this section, we outline how to use Theorem \ref{count} to deduce Theorem \ref{main}, following Shkredov's argument in \cite{Sh26}. \par

The starting point is that, a reduced fraction $\frac{a}{q}<1$ having bounded partial quotients is equivalent to $\frac{a}{q}$ satisfying certain Diophantine condition. Specifically, if all partial quotients of $\frac{a}{q}$ are bounded by $M$, then for all $1\leq x< q$, we have $x|ax|_q>\frac{q}{4M}$, where $|ax|_q$ is the distance of $ax$ to a nearest integer in $q\mathbb Z$. Conversely, if for all $1\leq x< q$, we have $x|ax|_q>\frac{q}{M}$, then all partial quotients of $\frac{a}{q}$ are bounded by $M$ (see Lemma 2, \cite{MMS26} for a proof).  \par

Given $\frac{e}{f}=[c_1, c_2, \cdots, c_s]<1$, for each $1\leq \nu\leq s$, we let $\frac{e_\nu}{f_\nu}=[c_1, \cdots, c_\nu] $ be the $\nu$-th convergent of $\frac{e}{q}$. The denominator $f_\nu$ is also called $\nu$-th continuant of $\frac{e}{f}$. 

Now for a parameter $t$, we consider the set $Q_M(t)$ consisting of all $\frac{e}{f}$ such that all partial quotients of $\frac{e}{f}$ up to $\nu$ are bounded by $M$, where $\nu$ is taken so that $f_\nu$ is the largest continuant of $\frac{e}{f}$ smaller than $t$. The set $Q_M(t)$ has a fractal structure, consisting of $\asymp t^{2w_M}$ many intervals of length $\asymp t^{-2}$. Now for a fixed $q$, we want each interval to contain many rationals with denominator $q$, so we put on the restriction $t<q^{\frac{1}{2}-\epsilon_0}$, and we let $Z_M(t)$ be the collection of numerators $a$ of all $\frac{a}{q}\in Q_M(t)$. \par

Inheriting from the fractal structure of $Q_M(t)$, the set $Z_M(t)$ contains a set of the form $K=\Lambda_1 \dotplus [1,N]$, where $|\Lambda_1|\asymp t^{2w_M}$ and $t^2 N\asymp q$. \par
If $a\in Z_M(t)$, it implies that for all $1<x\leq t$, $x|ax|_q>\frac{q}{M}$. If $a^{-1}(\mod q)\in Z_M(t)$ as well, then the inequality $x|ax|_q>\frac{q}{M}$ also holds for $x\in[\frac{q}{4Mt}, q)$. So we are looking for $a\in K\cap K^{-1}$, which naturally leads to the consideration of equation \eqref{0102}, for which Theorem \ref{count} gives an asymptotic counting for the solutions. In the paper \cite{MMS26}, the parameters $\epsilon_0, M$ are chosen in a way so that  $\frac{q}{4Mt}<t$ and the error term in \eqref{0905} is dominated by the main term, which leads to $M=O(\log q/\log\log q)$ and $\epsilon_0 \asymp \frac{1}{M}$. \par 
If one wants $M$ to be an absolute number, which requires $\epsilon_0\gg 1$, inevitably one needs to analyse the range $x\in [t, \frac{q}{4Mt}]$, and this is where the main contribution of \cite{Sh26} lies. For this, one needs to introduce a new parameter $H=N^{\frac{9}{20}}$, and consider a refined set $A=Z_M(t)\cap Z_{M^*}(tH)$ for some $M^*\approx 10 M$, which contains a similar fractal structure $\Lambda_2\dotplus [1, \frac{N}{H^2}]=\sqcup_{\lambda\in \Lambda_2}I_\lambda$, where $I_\lambda=\lambda+[0, N^{\frac{1}{10}}]$. \par
Take $N^*=N^{\frac{1}{100}}$. An application of Theorem \ref{count} shows that most $I_\lambda$ are $N^*$-good, in the sense that if dividing $I_\lambda$ into subintervals of length $N^*$, at least $(1-(N^*)^{-\frac{\eta}{2}})$ proportion of all such subintervals intersect with $A^{-1}$. Take any such $N^*$-good interval $I$. Let $I_1$ be a general subinterval of $I$ of length $N^*$ and let $a\in I_1$ such that $a^{-1}\in A$. Suppose all partial quotients of $\frac{a}{q}$ are bounded by $\tilde M$ for some absolute $\tilde M\geq \max\{M, M_*\}$, then $a$ is a desired numerator and we are done. Otherwise, $\frac{a}{q}$ has a partial quotient $c_\nu>\tilde M$ with corresponding continuant $x=x(a)$. Then we must have $x(a)\in[t, \frac{q}{4tM}]$. The denominator $x(a)$ is also called an \emph{$\tilde M$-critical denominator} of $\frac{a}{q}$. Furthermore, if all such $I_1$ contains a critical denominator, it is then possible to find $a, b, c\in I$, well separated in the sense that $b-a, c-b=\Theta(N^{\frac{1}{10}})$, with their $\tilde M$-critical denominators $x, y, z$. Since $a, b, c\in A\cap A^{-1}$, we can assume without loss of generality that $y, z\leq \sqrt{q}$. On the one hand, critical denominators exhibit some repelling phenomenon: There exists $x'\in [tN^{\frac{11}{20}}, tN^{\frac{3}{4}}]$ which is a continuant of $\frac{a}{q}$, such that the equation 
\begin{align}\label{0719}
\alpha x'+\beta y+\gamma z=0 
\end{align}
does not admit a solution for $(\alpha, \beta, \gamma)$ with all $\alpha, \beta, \gamma$ small. On the other hand, we have $|by|_q\leq \frac{q}{\tilde My}$, and since $x$ and $z$ are close to $y$ (in the sense of Euclidean distance), $|bx|_q, |bz|_q$ are also small. It can then be shown that there exists small $m_1, m_2, m_3$ such that $0\leq |m_1bx'+m_2by+m_3bz|_q < t$. $|m_1bx'+m_2by+m_3bz|_q $ can not be zero; otherwise, \eqref{0719} has a small solution $(\alpha, \beta, \gamma)=(m_1, m_2, m_3)$. So $0< |m_1bx'+m_2by+m_3bz|_q < t$. Since $b^{-1}\in Z_M(t)$, we then obtain that $|m_1x'+m_2y+m_3z|_q>\frac{q}{4Mt}$, which leads to a contradiction, having in mind that $m_1, m_2, m_3$ are small. 

\par

\bibliographystyle{alpha}

\bibliography{boundedgenerationV3}

\end{document}